\documentclass[titlepage,11pt]{article}
\usepackage{latexsym}
\usepackage{amsfonts}
\usepackage{amsmath,enumerate}
\usepackage{hyperref, fullpage, setspace}
\usepackage{amsthm, amssymb}
\usepackage{graphics}
\usepackage{graphicx}
\usepackage{epstopdf}

\newtheorem {theorem}            {}[section]

\newcommand\blackslug{\hbox{\hskip 1pt \vrule width 4pt height 8pt depth 1.5pt
        \hskip 1pt}}
\newcommand\bbox{\hfill \quad \blackslug \bigbreak}

\title{The Minimal Automorphism-Free Tree}
\author{
Ilhee Kim \qquad Ringi Kim \qquad Paul Seymour\\
\\
Princeton University\\
Princeton, NJ, USA
}

\date{January 20, 2013; revised \today}

\begin{document}
\maketitle

\begin{abstract}
A finite tree $T$ with $|V(T)| \geq 2$ is called {\it automorphism-free} if
there is no non-trivial automorphism of $T$.
Let $\mathcal{AFT}$ be the poset with the element set of all finite automorphism-free trees (up to graph isomorphism)
ordered by $T_1 \preceq T_2$ if $T_1$ can be obtained from $T_2$
by successively deleting one leaf at a time in such a way that each intermediate tree is also automorphism-free.
In this paper, we prove that $\mathcal{AFT}$ has a unique minimal element.
This result gives an affirmative answer to the question asked by Rupinski in \cite{Rup}.
\end{abstract}

\section{Introduction}
In this paper, every graph is finite and simple.
For a graph $G$, a bijection $\phi : V(G) \rightarrow V(G)$ is an {\it automorphism} of $G$
if $uv \in E(G)$ if and only if $\phi(u)\phi(v) \in E(G)$.
For example, the identity function on $V(G)$ is an automorphism of $G$. We call this identity function
the {\it trivial} automorphism.

For a tree $T$ with $|V(T)| \geq 2$, we say $T$ is {\it automorphism-free} if there is no non-trivial automorphism of $T$.
For instance,
the following graph $E_7$ in Figure 1 is automorphism-free.
It is easy to check that there are no automorphism-free trees with fewer than 7 vertices.
\newline

\begin{figure}[h!]
  \begin{center}
    \scalebox{0.8}
      {
        \includegraphics{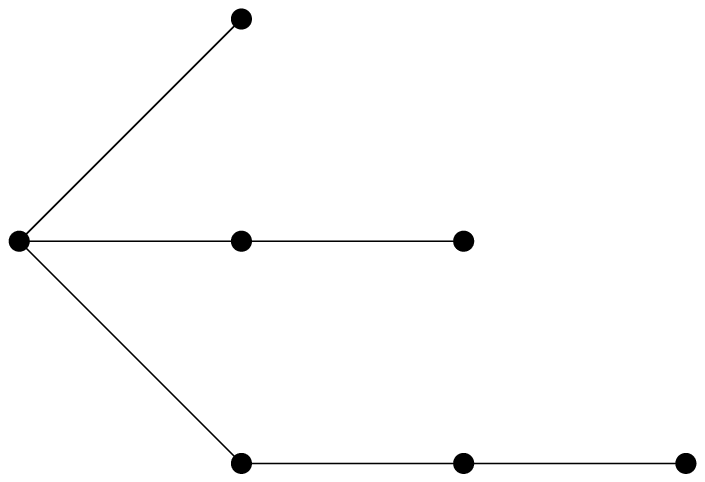}
      }
  \end{center}
  \caption{$E_7$}
\end{figure}

Let $\mathcal{AFT}$ be the poset (partially ordered set) with the element set of all finite automorphism-free trees
(up to graph isomorphism) ordered by $T_1 \preceq T_2$ if $T_1$ can be obtained from $T_2$ by successively deleting one leaf at a time
in such a way that each intermediate tree is also automorphism-free.

In this paper, we prove that $\mathcal{AFT}$ has a unique minimal element, namely $E_7$.

\begin{theorem}\label{main1}
Let $T$ be a minimal element of the poset $\mathcal{AFT}$. Then $T$ is isomorphic to $E_7$.
\end{theorem}

Equivalently,

\begin{theorem}\label{main2}
Every automorphism-free tree $T$ can be obtained  from $E_7$ by successively adjoining a leaf at a time
in such a way that each intermediate tree is also automorphism-free.
\end{theorem}

or,

\begin{theorem}\label{main3}
For every  automorphism-free tree $T$, $E_7$ can be obtained from $T$ by successively deleting a leaf at a time
in such a way that each intermediate tree is also automorphism-free.
\end{theorem}

This result gives an affirmative answer to the question asked by Rupinski in \cite{Rup}.
First, we start with some definitions.
A {\it component} of a graph $G$ is a maximal non-null subgraph of $G$.
For a vertex $u$ of a graph $G$, $G \setminus u$ denotes the graph obtained from $G$ by deleting the vertex $u$ (deleting all the edges incident with $u$ as well).
For an edge $uv$ of a graph $G$, $G \setminus uv$ denotes the graph obtained from $G$ by deleting the edge $uv$ (not
deleting the vertex $u$ or $v$).
For a vertex set $S \subseteq V(G)$ of a graph $G$, $G|S$ denotes the subgraph of $G$ induced by $S$.
For a tree $T$, a {\it leaf} $l$ is a vertex of degree one in $T$, and $p(l)$ denotes the (unique) neighbor of $l$ in $T$.
For a path $P$, the {\it length} of $P$ is the number of edges in $P$.
For a tree $T$ and $u,v \in V(T)$, $dist_T(u,v)$ is the length of the (unique) path from $u$ to $v$ in $T$.
For each $v \in V(T)$,  $d_T(v) := \max_{u \in V(T) \setminus \{v\}} dist_T(u,v)$.
We say $v \in V(T)$ is a {\it center} of $T$ if $d_T(v) \leq d_T(u)$ for every $u \in V(T)$,
and the {\it radius} $r(T)$ of $T$ denotes the number $d_T(v)$.

For the proof of \ref{main1}, we look at a minimal element $T$ of $\mathcal{AFT}$.
In $T$, we choose special leaves $l_1$ and $l_2$ by certain methods, and use the fact that
both $T \setminus l_1$ and $T \setminus l_2$ have non-trivial automorphisms.
From this, we find various properties that $T$ must have.
For instance, we prove that $T$ must have two centers, and $T \setminus l_1$ must have exactly one center,
and $T \setminus l_2$ must have two centers, etc.
Eventually we prove that $T$ must be isomorphic to $E_7$.
\newline

\section{Main proof}
\bigskip

The following is an easy lemma about centers in a tree. We omit the proof.
\bigskip

\begin{theorem}\label{centers}
Let $T$ be a tree with $|V(T)| \geq 2$. Let $l$ be a leaf of $T$ and let $\phi$ be an automorphism of $T$.
\begin{itemize}
\item[(1)] If $u$ and $v$ are distinct centers of $T$, then $uv \in E(T)$. In particular, there are at most two centers of $T$.
\item[(2)] If $u$ and $v$ are distinct centers of $T$, then every path of length $r(T)$ from $u$ contains $v$, and vice versa.
\item[(3)] If $u$ is the unique center of $T$, then $\phi(u) = u$.
\item[(4)] If $u$ and $v$ are distinct centers of $T$, then $\phi$ either fixes $u$ and $v$ or switches them.
\item[(5)] If $u$ is the unique center of $T$, then it is a center of $T \setminus l$ as well.
\item[(6)] If $u$ and $v$ are distinct centers of $T$, then every center of $T \setminus l$ is either $u$ or $v$.
\end{itemize}
\end{theorem}

\bigskip

For a given tree $T$ with $|V(T)| \geq 2$ and a vertex $u$ of $T$,
we say a leaf $l (\neq u)$ is
a {\it special leaf} with respect to $T$ and $u$ if the following statement holds for $l$.
\bigskip

Let $P$ be the path from $u$ to $l$ in $T$ and number the vertices of $P$ as
$v_1 (=u), \cdots, v_m (=l)$ in order. For each $i = 1,\cdots,m-1$,
let $C_i$ be the component of $T \setminus v_i$ containing $l$.
Then for every component $C$ of $T \setminus v_i$ not containing $u$, $|V(C)| \geq |V(C_i)|$.
\bigskip

\begin{theorem}\label{special}
Let $T$ be a tree with $|V(T)| \geq 2$ and let $u \in V(T)$.
Then, there exists a special leaf with respect to $T$ and $u$.
\end{theorem}

{\noindent \bf Proof.}
We proceed by induction on $|V(T)|$. It is easy to see that the statement holds for $|V(T)| = 2$.
Consider all neighbors of $u$. Each one is in its own component of $T \setminus u$.
Among those components, we take one with the least number of vertices.
(If there is more than one smallest component, just pick any one of them.)
Let $C$ be the component of $T \setminus u$ we chose and let $v$ be the neighbor of $u$ in $C$.
Now, look at all children of $v$ (the neighbors of $v$ in $C$). If there are no children of $v$,
then $v$ is a special leaf with respect to $T$ and $u$ we are looking for.
Therefore we may assume $|V(C)| \geq 2$. Then from the induction hypothesis, there is a special leaf $l$ with respect to $C$ and $v$.
Then, it is easy to check that $l$ is a special leaf with respect to $T$ and $u$ as well.
This proves \ref{special}. \bbox

\begin{theorem}\label{AFT1}
Let $T$ be a minimal tree in the poset $\mathcal{AFT}$. Then $T$ has two centers.
\end{theorem}

{\noindent \bf Proof.} For the sake of contradiction, suppose $T$ has only one center $u$.
Let $l_1$ be the special leaf with respect to $T$ and $u$. Let $T' = T \setminus l_1$ and take
a non-trivial automorphism $\phi$ of $T'$. By \ref{centers} (5), $u$ is a center of $T'$ as well.
\newline

\noindent (1) {\it  $\phi$ does not fix $u$. In particular, there is another center of $T'$.}
\newline
\newline
Suppose $\phi$ fixes $u$.
Notice that $\phi$ does not fix $p(l_1)$ because otherwise we can extend $\phi$ to a non-trivial automorphism of $T$
by assigning $\phi(l_1) = l_1$. In particular, $u \neq p(l_1)$.
Let $P$ be the path from $u$ to $p(l_1)$ in $T'$, and number the vertices of $P$ as $v_1 (=u), \cdots, v_m (= p(l_1))$ in order ($m \geq 2$).
Let $k$ be the largest number such that $\phi(v_k) = v_k$ (and so, $k \leq m-1$ and $\phi(v_{k+1}) \neq v_{k+1}$).

In $T' \setminus v_k$, let $C_1$ be the component containing $v_{k+1}$ (this component contains $p(l_1)$ as well),
and let $C_2$ be the component containing $\phi(v_{k+1})$.
Clearly, $C_1$ and $C_2$ are different since $v_{k+1}$ and $\phi(v_{k+1})$ are both neighbors of $v_k$.

\begin{figure}[h!]
  \begin{center}
    \scalebox{1.0}
      {
        \includegraphics{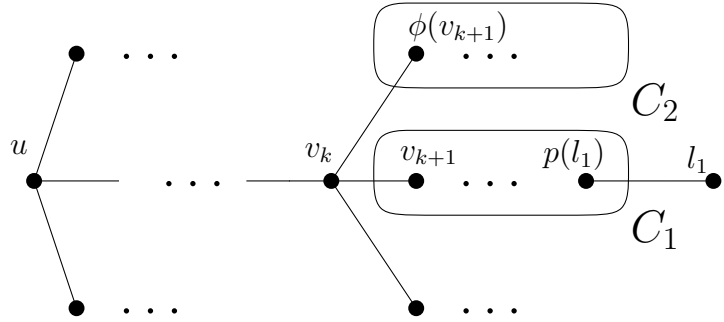}
      }
  \end{center}
  \caption{$C_1$ and $C_2$}
\end{figure}

Notice that $\phi | V(T' \setminus v_k)$ is an automorphism of $T' \setminus v_k$ since $\phi$ fixes $v_k$.
This implies $C_1$ and $C_2$ are isomorphic. In particular, $|V(C_1)| = |V(C_2)|$.
But back in $T$, $T|(V(C_1) \cup \{l_1\})$ and $C_2$ are two components of $T \setminus v_k$, and $|V(C_1) \cup \{l_1\}| > |V(C_2)|$. This contradicts the definition of $l_1$. This proves (1).
\newline

Let $v$ be the center of $T'$ different from $u$.
\newline

\noindent (2) {\it $d_T(v) = r(T) + 1$, and there is a unique path of length $d_T(v)$ from $v$ in $T$, namely the path from $v$ to $l_1$.}
\newline
\newline
In $T$, $u$ is the unique center. But in $T'$, both $u$ and $v$ are centers.
Therefore
$$d_T(v) - 1 \geq r(T) = d_T(u) \geq d_{T'}(u) = d_{T'}(v) \geq d_T(v) - 1.$$
In particular, $d_T(v) - 1 = d_{T'}(v)$. This implies there is no path of length $d_T(v)$ from $v$ in $T'$.
Since there is a path of length $d_T(v)$ from $v$ in $T$, namely the path from $v$ to $l_1$, it must be unique.
This proves (2).
\newline

\noindent (3) {\it $l_1$ is not a neighbor of $u$.}
\newline
\newline
Suppose $l_1$ is a neighbor of $u$.
Then from (2),
$$r(T') = d_{T'}(v) = d_T(v) - 1 = dist_T(v, l_1) - 1 = 1.$$
Since $T'$ has two centers and $r(T') = 1$, $|V(T')| = 2$ and $|V(T)| = 3$.
But this is impossible since there is no automorphism-free tree with three vertices.
This proves (3).
\newline

 \noindent (4) {\it $p(l_1)$ has degree two in $T$.}
\newline
\newline
If there exists another child $w$ of $p(l_1)$ in $T$,
then the path from $v$ to $w$ is another path of length $d_T(v)$ from $v$ in $T$, which is impossible by (2).
Therefore $l_1$ is the unique child of $p(l_1)$ in $T$.
This proves (4).
\newline

Let $T_u$ and $T_v$ be the two components of $T \setminus uv$
containing $u$ and $v$, respectively.
Note that $\phi$ switches $u$ and $v$ by (1).
And $T_u \setminus l_1$ is isomorphic to $T_v$ by $\phi$.
Let $l_2$ be a special leaf with respect to $T_v$ and $v$
($l_2$ exists since $V(T_v) \geq 2$).
Let $T'' = T \setminus l_2$ and take a non-trivial automorphism $\psi$ of $T''$.
\newline

\begin{figure}[h!]
  \begin{center}
    \scalebox{0.6}
      {
        \includegraphics{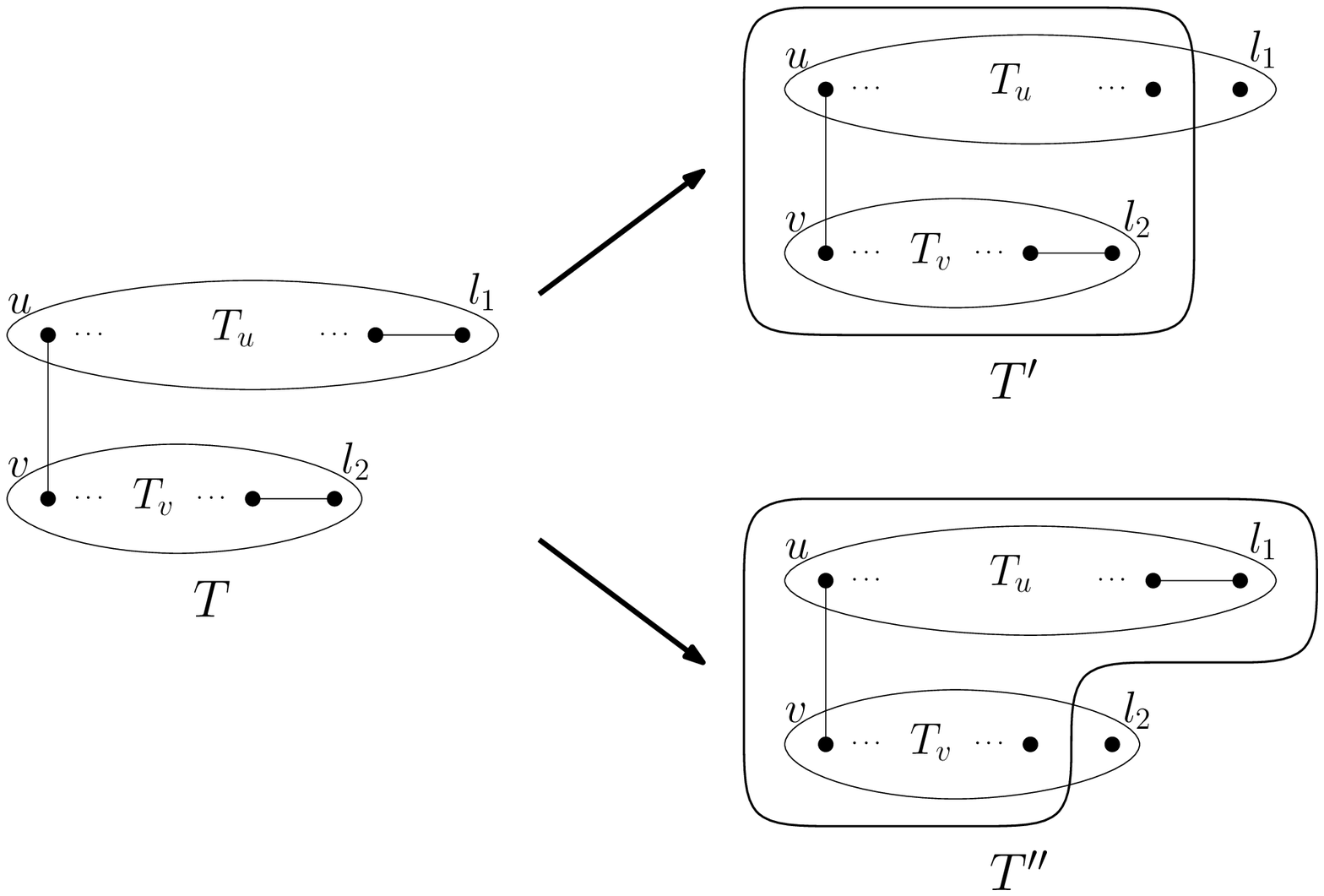}
      }
  \end{center}
  \caption{$T'$ and $T''$}
\end{figure}

\noindent (5) {\it $u$ is a center of $T''$, and $v$ is not.}
\newline
\newline
Again, $u$ is a center of $T''$ by \ref{centers} (5).
But, $v$ is not a center of $T''$ because from (2), $$d_{T''}(v) \geq dist_{T''}(v, l_1) = r(T) + 1 > d_T(u) \geq d_{T''}(u).$$
This proves (5).
\newline

\noindent (6) {\it $\psi$ does not fix $v$. Moreover, $\psi(V(T_v \setminus l_2)) \subseteq V(T_u \setminus u)$.}
\newline
\newline
Suppose $\psi$ fixes $v$. Then $u$ is fixed as well because
among the neighbors of $v$, $u$ is the unique center of $T''$ (although $u$ might not be the unique center of $T''$).
Since both $u$ and $v$ are fixed by $\psi$, $\psi(V(T_v)) = V(T_v)$.
Therefore
$\psi | V(T_v \setminus l_2)$ is a non-trivial automorphism of $T_v \setminus l_2$
(otherwise we can extend $\psi$ to a non-trivial automorphism of $T$ by assigning $\psi(l_2) = l_2$).
Then by the same argument as in (1),
this contradicts the definition of $l_2$.

Since $v$ is adjacent to a center of $T''$, so is $\psi(v)$.
And $v$ is the unique such vertex in $V(T_v \setminus l_2)$.
Therefore $\psi(v)$ does not belong to $V(T_v \setminus l_2)$ because $\psi(v) \neq v$.
Also, every member of $\psi(V(T_v \setminus l_2))$ does not belong to $V(T_v \setminus l_2)$ either
because $T_v \setminus l_2$ is a component of $T'' \setminus u$ containing $v$.
Therefore $\psi(V(T_v \setminus l_2)) \subseteq V(T_u \setminus u)$.
This proves (6).
\newline

Let $C$ be the component of $T'' \setminus u$ containing $\psi(V(T_v \setminus l_2))$.
Let $n = |V(T_v)|$.
\newline

\noindent (7) {\it $|V(C)|$ is either $n$ or $n-1$.}
\newline
\newline
Since $\psi(V(T_v \setminus l_2)) \subseteq V(C)$, $|V(C)| \geq n-1$ .
Recall that $|V(T_u)| - 1 = |V(T_u \setminus l_1)| = |V(T_v)| = n$ and $V(C)$ is a subset of $V(T_u \setminus u)$.
Therefore
$$|V(C)| \leq |V(T_u \setminus u)| = |V(T_u)| - 1 = n.$$
This proves (7).
\newline

 \noindent (8) {\it $|V(C)| = |V(T_u \setminus u)|= n$. In particular, the degree of $u$ is two (both in $T''$ and $T$), and $T''\setminus u$ consists of two components $T_u \setminus u (= C)$ and $T_v \setminus l_2$. }
\newline
\newline
From (3), $u$ has no neighbor which is a leaf (both in $T$ and $T''$).
In particular, every component of $T'' \setminus u$ has more than one vertex.
Note that the union of all components of $T''\setminus u$ different from $T_v\setminus l_2$ has size $|V(T_u \setminus u)| = n$.
Since $C$ is one of them whose size is at least $n-1$, there cannot be another one (and so, $|V(C)| = n$).
Therefore $u$ has degree two, and this proves (8).
\newline

Notice that the degree of $v$ in $T$ is also two since $v = \phi(u)$.
\newline

\noindent(9) {\it $\psi$ does not fix $u$. In particular, $T''$  has two centers $u$ and $\psi(u)$, and  $T_v\setminus l_2 \cong T_u \setminus u \setminus \psi(u)$.}
\newline
\newline
By (5), $v$ is not a center of $T''$ and it is adjacent to a center of $T''$.
Therefore $\psi(v)$ has the same property in $T''$.
But $\psi(v)$ is not adjacent to the center $u$, because if it is, then $T_u \setminus u$ is isomorphic to $T_v \setminus l_2$
and hence, $|V(C)| = |V(T_v \setminus l_2)| = n-1$. This is impossible by (8).

Therefore $\psi(v)$ is adjacent to another center, and this also implies $\psi(u)\neq u$.
Since two centers are adjacent, $\psi(u)$ must be the neighbor of $u$ different from $v$.
Since $\psi(u)$ also has degree two, the neighbor of $\psi(u)$ different from $u$ is $\psi(v)$.
And $T_v\setminus l_2$ is isomorphic to $T_u \setminus u \setminus \psi(u)$. This proves (9).
\newline

\begin{figure}[h!]
  \begin{center}
    \scalebox{0.8}
      {
        \includegraphics{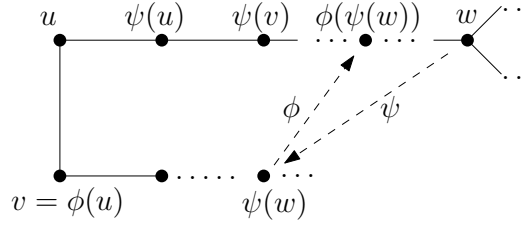}
      }
  \end{center}
  \caption{$\phi$ and $\psi$}
\end{figure}


Note that 
$T_u \setminus l_1 \cong T_v$ by $\phi$, and
$T_u \setminus u \setminus \psi(u) \cong T_v \setminus l_2$ by $\psi$.
\newline

\noindent (10) {\it $T_u$ and $T_v$ are paths.}
\newline
\newline
It is enough to show that $T_u$ is a path since $T_v \cong T_u \setminus l_1$.
Suppose there exists a vertex of degree at least three in $T_u$.
Choose such a vertex $w \in V(T_u)$ with $dist_{T}(u,w)$ as small as possible.
Then, $\phi(\psi(w))$ has degree at least three in $T_u$ as well. To see this,
first observe that $\psi(w) \in V(T_v) \setminus \{v,l_2\}$ has degree at least three since 
$w \in V(T_u) \setminus \{u, \psi(u), \psi(v)\}$, and $T_u \setminus u \setminus \psi(u) \cong T_v \setminus l_2$.
Therefore $\phi(\psi(w)) \in V(T_u)$ has degree at least three since $T_v \cong T_u \setminus l_1$. 
But then,
$$dist_{T}(u,\phi(\psi(w)))=dist_T(\phi^{-1}(u),\psi(w))=dist_{T}(v,\psi(w))=dist_{T}(\psi^{-1}(v),w)<dist_{T}(u,w).$$
This contradicts our choice of $w$.
This proves (10).
\newline



Since both $u$ and $v$ have degree two in $T$, $T$ is a path by (10).
This contradicts the fact that $T$ is automorphism-free.
Therefore $T$ has two centers.
This proves \ref{AFT1}. \bbox
\bigskip

{\noindent \bf Proof of \ref{main1}.}  By \ref{AFT1}, $T$ has two centers $u$ and $v$.
Let $T_u$ and $T_v$ be the two components of $T \setminus uv$ containing $u$ and $v$, respectively.

Let $l_1$ be a special leaf with respect to $T_u$ and $u$ and let $l_2$ be a special leaf with respect to $T_v$ and $v$.
Let $x$ be the shortest distance from $u$ to a vertex in $T_u$ whose degree in $T$ is at least three.
If there is no such vertex, then set $x$ as $\infty$.
Similarly, let $y$ be the shortest distance from $v$ to a vertex in $T_v$ whose degree in $T$ is at least three.

Without loss of generality, we may assume $|V(T_u)| \geq |V(T_v)|$. And further we may assume
if $|V(T_u)| = |V(T_v)|$ then $x \geq y$ by switching $u$ and $v$ if necessary.
We first consider $T' = T \setminus l_1$.
Let $\phi$ be a non-trivial automorphism of $T'$.
\newline

\noindent (1) {\it $u$ and $v$ are centers of $T'$.}
\newline
\newline
By \ref{centers} (2), $d_{T'}(u) = d_T(u)$ since every path of length $d_T(u)$ from $u$ in $T$ does not contain $l_1$.
And by \ref{centers} (6), $v$ is a center of $T'$ since $d_{T'}(v) \leq d_{T}(v) = d_T(u) = d_{T'}(u)$.
Again by \ref{centers} (6), if there is a center of $T'$ different from $v$, then it must be $u$.
For the sake of contradiction, suppose $v$ is the unique center of $T'$.
Then $p(l_1)$ is a leaf in $T'$ by the same argument as in (4) in the proof of \ref{AFT1}.
Then, $\phi$ does not fix $u$ since otherwise it contradicts the definition of $l_1$, by the same argument as in (6) in the proof of \ref{AFT1}.
Therefore in $T' \setminus v$, the component $T_u \setminus l_1$ is isomorphic to another component $C$ of $T' \setminus v$.
Note that $$|V(C)| = |V(T_u \setminus l_1)| = |V(T_u)| - 1.$$
Since $V(C)$ is a subset of $V(T_v \setminus v)$,
$$|V(T_v)| \geq |V(C)| + 1 = |V(T_u)|.$$
Together with our assumption $|V(T_u)| \geq |V(T_v)|$, $|V(T_v)| = |V(T_u)|$. Moreover $T' \setminus v$ has exactly two components, namely $C$ and $T_u \setminus l_1$.
In particular, $v$ has degree two in $T$.

Next, there exists a vertex in $V(T_v)$ whose degree in $T$ is at least three
because otherwise $T_v$ is a path, and hence so is $T'$, and so is $T$ because $p(l_1)$ has degree two in $T$.
But then, $y = x+1$ by $\phi$, and this contradicts our assumption that $x \geq y$ if $|V(T_u)| = |V(T_v)|$.
Therefore $v$ is not the unique center of $T'$.
This implies $u$ is a center of $T'$ as well.
This proves (1).
\newline

\noindent (2)  {\it $\phi$ switches $u$ and $v$, and $|V(T_v)| = |V(T_u)| -1$.}
\newline
\newline
Again, if $\phi$ fixes $u$, then $\phi$ fixes $v$ as well and this contradicts our choice of $l_1$.
Since $\phi$ switches $u$ and $v$, $T_v$ and $T_u \setminus l_1$ are isomorphic.
In particular, $|V(T_v)| = |V(T_u)| - 1$.
This proves (2).
\newline

Now we consider $T'' = T \setminus l_2$.
Let $\psi$ be a non-trivial automorphism of $T''$.
\newline

\noindent (3) {\it $\psi$ does not fix $v$, and $u$ is the unique center of $T''$.}
\newline
\newline
By the same argument as in the proof of (1),
$u$ is a center of $T''$, and if there is another one, then it must be $v$.
Suppose $\psi$ fixes $v$. Then, again $u$ is fixed as well and this contradicts the definition of $l_2$
by the same argument as in (6) in the proof of \ref{AFT1}.
Therefore $\psi$ does not fix $v$.

Now, suppose $v$ is a center of $T''$.
Since $\psi$ does not fix $v$, it switches $u$ and $v$ and $T_v \setminus l_2$ is isomorphic to $T_u$.
But then, $|V(T_v \setminus l_2)| =  |V(T_v)| - 1 = |V(T_u)| - 2 \neq |V(T_u)|$,
a contradiction. This proves (3).
\newline
\begin{figure}[h!]
  \begin{center}
    \scalebox{0.8}
      {
        \includegraphics{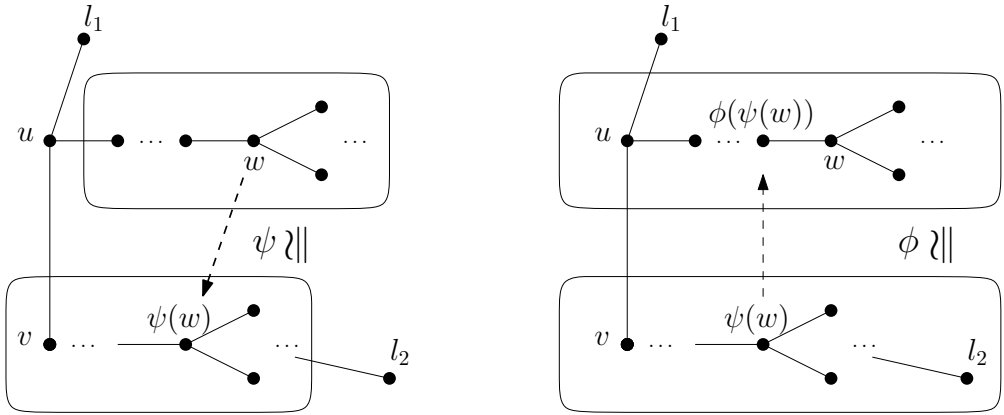}
      }
  \end{center}
  \caption{$\phi$ and $\psi$}
\end{figure}

Since $u$ is the unique center of $T''$ and $\psi$ does not fix $v$, the component $T_v \setminus l_2$ of $T'' \setminus u$ is isomorphic to another component $C$ of $T'' \setminus u$.
Note that the union of all components of $T''\setminus u$ different from $T_v \setminus l_2$ is exactly $T_u \setminus u$.
And $C$ has size $|V(T_v)| - 1 = |V(T_u)| - 2$.
This implies that there are exactly three components of $T''\setminus u$, namely $T_v \setminus l_2$, $C$, and the third one with a single vertex.
Therefore $u$ has a neighbor of degree one, and this implies that $l_1$ is a neighbor of $u$.
\newline

Now, $T_v \cong T_u \setminus l_1$ by $\phi$, and  $T_u \setminus u \setminus l_1 \cong T_v \setminus l_2$ by $\psi$.
\newline

\noindent (4) $T_u  \setminus l_1 $ and $T_v$ are paths.
\newline
\newline
It is enough to show that $T_u  \setminus l_1$ is a path since $T_v  \cong T_u \setminus l_1$ by $\phi$.
Suppose there exists a vertex of degree at least three in $T_u  \setminus l_1$. Choose such a vertex $w \in V(T_u  \setminus l_1)$
with $dist_T(u,w)$ as small as possible.
Then, $\phi(\psi(w))$ has degree at least three in $T_u\setminus l_1$ as well.
To see this, first observe that $\psi(w) \in V(T_v \setminus l_2)$ has degree at least three in $T_v$ since 
$w \in V(T_u\setminus u)$ and $T_u \setminus u \setminus l_1 \cong T_v \setminus l_2$.
Therefore $\phi(\psi(w)) \in V(T_u \setminus l_1)$ has degree at least three in $T_u \setminus l_1$ since
$T_u \setminus l_1 \cong T_v$.
But then,
$$dist_T(u,\phi(\psi(w)))=dist_T(\phi^{-1}(u),\psi(w)) = dist_T(v,\psi(w))$$
$$=dist_T(\psi^{-1}(v),w)=dist_T(\psi(v),w)<dist_T(u,w).$$
This contradicts our choice of $w$.
This proves (4).
\newline

By (4), $T$ is a tree with a unique vertex of degree three, namely $u$, and one of the three components of $T \setminus u$ consists of a single vertex, namely $l_1$, and the other two components $T_v$ and $T_u \setminus u \setminus l_1$ are paths.
Let $|V(T_u \setminus u \setminus l_1)| = k$; then $|V(T_v)| = k+1$ by (2).
Finally if $k > 2$, then deleting the leaf of $T$ in $V(T\setminus u \setminus l_1)$ yields another automorphism-free tree, and if $k = 1$, then $|V(T)| = 5$, and so $T$ is not automorphism-free.
Therefore $k = 2$, and $T$ is isomorphic to $E_7$.
This proves \ref{main1}. \bbox



\end{document}